

\magnification=\magstep1 \baselineskip=15pt
\parskip=4pt

\def\<{\langle}
\def\>{\rangle}
\def\z{\zeta}
\def\f{\phi}

\def\E{{\cal E}}
\def\F{{\cal F}}
\def\G{{\cal G}}

\def\C{{\bf C}}
\def\R{{\bf R}}
\def\d{\partial}
\def\Re{{\rm Re}}

\def\Conv{{\rm Conv}}

\font\auth=cmcsc10
\def\jour{\sl}
\def\vol{\bf}
\def\bysame{\vrule height0.4pt depth0pt width1truein \ , }

\topglue 1truein
\centerline{\bf
EXTREMAL DISCS AND ANALYTIC CONTINUATION}
\centerline{\bf
OF PRODUCT CR MAPS}
\bigskip
\centerline{{\auth A.~Scalari and A.~Tumanov}}
\bigskip

\beginsection  Introduction

One of the essentially multidimensional phenomena in complex
analysis is the forced analytic continuation of a germ of
a biholomorphic map $M_1\to M_2$ between real analytic
manifolds $M_1$ and $M_2$ in $\C^n$, $n>1$.
Poincar´e (1907) observed that a biholomorphic map sending an
open piece of a unit sphere in $\C^2$ to another such open
piece must be an automorphism of the unit ball. This
was proved for $\C^n$ by Tanaka (1962) and then rediscovered
by Alexander [A].

Pinchuk [P] proved that if $M_1$ and $M_2$ are strictly pseudoconvex
real analytic nonspherical hypersurfaces and $M_2$ is compact,
then a germ of a biholomorphic map $M_1\to M_2$ holomorphically
extends along any path in $M_1$. Ezhov, Kruzhilin and Vitushkin [EKV]
gave a different proof of that result.
Webster [W] proved that a germ of a biholomorphic map
$M_1\to M_2$ between real algebraic Levi non-degenerate
hypersurfaces in $\C^n$ is algebraic.

There is an impressive number of publications
in which $M_1$ and $M_2$ are real
{\it algebraic} manifolds of different dimensions or higher
codimension, in particular real quadratic
manifolds (see [BER]).
Hill and Shafikov [HS] prove the analytic continuation result
in higher codimension in which only one of the manifolds
$M_1$ and $M_2$ is assumed to be algebraic.
There are many more results on the problem that we
omit here, see e.g. [BER], [HS] for references.

Despite the large amount of work done on the problem,
seemingly there are no results in the literature in which
$M_1$ and $M_2$ are manifolds of higher codimension in $\C^n$
and neither of them is algebraic.
In this paper we consider the case in which
$M_1$ is a real analytic strictly pseudoconvex manifold
and
$M_2$ is the cartesian product of two or more compact
strictly convex real analytic hypersurfaces.
For the case in which $M_2$ is the product of two spheres,
the result was obtained earlier by the first author [Sc].
In this paper we significantly simplify and generalize
the proof given in [Sc]. Following [Sc] we use a new method
based on extremal discs in higher codimension.
As byproducts, we obtain some properties of extremal discs
that may be used elsewhere.

\beginsection  1. Strictly pseudoconvex manifolds

In this section we recall basic notations and definitions
concerning real manifolds in complex space.

Let $M$ be a $C^\infty$ smooth real generic
manifold in $\C^N$ of real codimension $k$.
Recall that $M$ is {\it generic} if
$T_p(M)+JT_p(M)=T_p(\C^N)$, $p\in M$, where
$T(M)$ denotes the tangent bundle to $M$, and
$J$ is the operator of multiplication by
the imaginary unit in $T(\C^N)$.
Recall the {\it complex tangent space} $T^c_p(M)$
of $M$ at $p\in M$ is defined as
$T^c_p(M)=T_p(M)\cap JT_p(M)$.
If $M$ is generic, then $M$ is a CR manifold,
which means that $\dim_\C T^c_p(M)$ is independent
of $p$, and $T^c(M)$ forms a bundle.
Recall the space
$T^{(1,0)}_p(M)\subset T_p(M)\otimes\C$ of complex
$(1,0)$-vectors is defined as
$T^{(1,0)}_p(M)=
\{X\in T_p(M)\otimes\C: X=\sum a_j\,\d/\d z_j \}$.
The CR dimension $\dim_{CR}(M)$ of $M$ is equal to
$\dim_\C T^c_p(M)=\dim_\C T^{(1,0)}_p(M)$.
If $\dim_{CR}(M)=n$, then $N=n+k$.

Let $T^*(\C^N)$ be the real cotangent bundle of $\C^N$.
Since every (1,0) form is uniquely determined by its
real part, we represent $T^*(\C^N)$ as the space
of (1,0) forms on $\C^N$.
Then $T^*(\C^N)$ is a complex manifold.
Let $N^*(M)\subset T^*(\C^N)$ be the real conormal bundle
of $M\subset\C^N$.
Using the representation of $T^*(\C^N)$ by (1,0) forms,
we define the fiber $N^*_p(M)$ at $p\in M$ as
$$
N^*_p(M)=\{\f\in T^*_p(\C^N):\Re\f|_{T_p(M)}=0\}.
$$

We use the angle brackets $\<,\>$ to denote the
natural pairing between vectors and covectors,
so we write $\<\f,\xi\>=\sum \f_j\xi_j$ for their
coordinate representations.

In a fixed coordinate system, we will identify
$\f=\sum \f_j\,dz_j\in T^*(\C^N)$ with the vector
$\f=(\f_1,\dots,\f_N)\in\C^N$.
Then for $\f\in N^*_p(M)$, the vector $\bar\f$ is
orthogonal to $M$ in the real sense, that is
$\Re\<\f,X\>=0$ for all $X\in T_p(M)$.

Since $M$ is generic, then locally $M$ can be
defined as
$\rho(z)=0$,
where
$\rho=(\rho_1,\dots,\rho_k)$
is a smooth real vector function such
that $\d\rho_1\wedge\dots\wedge\d\rho_k\ne 0$.
The forms
$\d\rho_j$, $(j=1,\dots,k)$,
define a basis of $N^*_p(M)$, so every
$\f\in N^*_p(M)$ can be written as
$\f=\sum c_j\d\rho_j$, $c_j\in\R$.

For every $\f\in N^*_p(M)$ we define the Levi form $L(p,\f)$ of $M$
at $p\in M$ in the conormal direction $\f=\sum c_j\partial\rho_j$ as
$$
L(p,\phi)(X,Y)=-\sum c_j \partial\bar\partial\rho_j(X,\bar Y),
$$
where $X,Y\in T^{1,0}_p(M)$. The form $L(p,\f)$ is a hermitian form
on $T^{1,0}_p(M)$. This definition is independent of the defining
function.
The forms $L(p,\f)$ can be regarded as components of the
$N_p(M)$-valued Levi form $L(p)$, where
$N(M)=T(\C^N)|_M/T(M)$ is the normal bundle of $M\subset\C^N$.
Indeed, $L(p)(X,X)\in N_p(M)$ is such an element that
$$
\Re\<\f, L(p)(X,X)\>=L(p,\f)(X,X) \quad
\hbox{for all} \quad
\f\in N^*_p(M).
$$
The Levi cone $\Gamma_p\subset N_p(M)$ is defined as the convex
span of the values of the Levi form $L(p)$, that is
$$
\Gamma_p=\Conv\{ L(p)(X,X): X\in T^{1,0}_p(M), X\ne0\}.
$$
We also need the Levi cone $H_p\subset T_p(M)$.
We put
$$
H_p=\{ \xi\in T_p(M): [J\xi]\in\Gamma_p \},
$$
where the brackets denote the class in the quotient space
$N_p(M)$.
If $M$ is a strictly pseudoconvex hypersurface, then
$\Gamma_p$ is the half-line defined by the inner normal
to $M$ at $p$ and $H_p$ is a half-space of $T_p(M)$.
The dual Levi cone $\Gamma_p^*$ is defined as
$$
\Gamma_p^*=\{ \phi\in  N^*_p(M): L(p,\phi)>0 \},
$$
where $L(p,\phi)>0$ means that the form $L(p,\phi)$
is positive definite.
The cones $\Gamma_p$ and $\Gamma^*_p$ are dual,
that is $\xi\in \Gamma_p$ iff $\Re\<\phi,\xi\>>0$
for all $\phi\in\Gamma^*_p$.

We say that $M$ is {\it strictly pseudoconvex} at $p$
if $\Gamma_p^*\ne \emptyset$.
We say that $M$ is strictly pseudoconvex if
it holds at every $p\in M$.
We say that the Levi form $L(p)$ is {\it generating} if
$\Gamma_p$ has nonempty interior.

Changing notation,
we introduce the coordinates $(z,w)\in\C^N$,
$z=x+iy\in\C^k$,
$w\in\C^n$, so that
the defining function of $M$ can be chosen in the form
$\rho=x-h(y,w)$,
where $h=(h_1,\dots,h_k)$ is a smooth real vector function,
and the equations of $M$ take the form (see, e.g., [BER])
$$
x_j=h_j(y,w)=\<A_jw,\bar w\>+O(|y|^3+|w|^3),
\qquad 1\le j \le k,
\eqno{(1.1)}
$$
where $A_j$ are hermitian matrices.
Then $T_0^{1,0}(M)$ is identified with the $w$-space $\C^n$ and for
$\f=\sum c_j dz_j\in N_0^*(M)$, the Levi form $L(0,\f)$ has the matrix
$\sum c_j A_j$.
Hence,
the manifold $M$ of the form (1.1) is strictly pseudoconvex at 0
if and only if there exists $c\in\R^k$ such that
$\sum c_j A_j>0$. It has a generating Levi form at 0 if
and only if the matrices $A_1,\dots,A_k$ are linearly
independent.

We say that a vector valued hermitian form $B$
{\it splits into scalar forms}
of dimensions $(n_1, ... , n_k)$
if the source and target spaces $V$ and $Z$ of $B$ split
into direct sums $V=\sum V_j$, $Z=\sum Z_j$, $\dim V_j=n_j>0$, $\dim Z_j=1$,
such that $B(u,v)=\sum B_j(u_j,v_j)$, where $u_j,v_j\in V_j$, $u=\sum u_j$,
$v=\sum v_j$, and $B_j$ is a $Z_j$-valued hermitian form on $V_j$.
We need the following simple
\medskip

{\bf Proposition 1.1.}
{\sl
Let $M$ be a connected
real analytic generic manifold in $\C^N$.
Suppose that the Levi form of $M$ splits into scalar forms
on an open subset of $M$. Then it splits into scalar forms
everywhere on $M$. If $M$ is strictly pseudoconvex,
then the Levi form is generating and splits into
positive definite forms.
}

{\it Proof.}
The set of all splittable hermitian forms is a real analytic
(even algebraic) subset of the set of all hermitian forms.
The map $M\ni p\mapsto L(p)$ is real analytic.
Since it takes an open set of $M$ to splittable forms
and since $M$ is connected, then the whole image belongs to
splittable forms. The rest of the conclusions hold
automatically. The proof is complete.

\beginsection  2. Extremal discs

We recall some facts of the theory of extremal discs [L], [T1].

Let $M$ be a smooth generic manifold in $\C^N$.
An {\it analytic disc} in $\C^N$ is a continuous
mapping $f:\bar\Delta\to\C^N$ holomorphic
in the unit disc $\Delta$.
We say that $f$ is {\it attached} to $M$
if $f(b\Delta)\subset M$.

An analytic disc $f$ attached to $M$ is
called {\it stationary} if there exists a
nonzero continuous holomorphic mapping
$f^*:\bar\Delta\setminus\{0\}\to T^*(\C^N)$,
such that $\tilde f=\z f^*$ is holomorphic in $\Delta$ and
$f^*(\z)\in N^*_{f(\z)}(M)$ for all $\z\in b\Delta$.
In other words, $f^*$ is a punctured analytic disc
with a pole of order at most one at zero attached to
$N^*(M)\subset T^*(\C^N)$ such that the
natural projection sends $f^*$ to $f$.
We call $f^*$ a {\it lift} of $f$, and we always use
the term ``lift'' in this sense.

We call a disc $f$ {\it defective}
if it has a nonzero lift
$f^*$ holomorphic in the whole
unit disc including 0.
For a strictly convex hypersurface, all
defective discs are constant.

We call a lift $f^*$ of a stationary disc $f$
{\it supporting} if for all
$\zeta\in b\Delta$, $f^*(\zeta)$ defines a (strong)
supporting real hyperplane to $M$ at $f(\zeta)$,
that is
$$
\Re \< f^*(\z), p-f(\z) \> \ge \epsilon |p-f(\z)|^2 \quad
\hbox{ for all }
\zeta\in b\Delta \hbox{ and }
p\in M,
\eqno(2.1)
$$
for some $\epsilon>0$. Stationary discs with supporting lifts have
important extremal properties, but we do not need them here.
Nevertheless, we call such $f$ {\it extremal} and we call
the pair $(f,f^*)$ an {\it extremal pair}.
Although $f$ is completely determined by $f^*$,
we prefer to use the excessive notation $(f,f^*)$, because it
lets us describe $f^*$ by its fiber coordinates in $T^*(\C^N)$.
Note that (2.1) implies that $f^*(\zeta)\in \Gamma^*_{f(\zeta)}$
for $\zeta\in b\Delta$.

If $M$ is the boundary of a strictly convex domain $D\subset\C^N$,
then the set of all extremal discs is smoothly parametrized by
the correspondence $f \leftrightarrow (f(0),f(1))\in D\times bD$.
The set of all extremal pairs is paramerized by
$D\times bD\times\R^+$ because the lift of an extremal
disc is unique up to a positive constant factor, see [L].

In higher codimension there is a local parametrization of
the set of extremal pairs.

{\bf Theorem 2.1} [T1]. {\sl
Let $M\subset\C^N$ be a smooth (resp. real analytic)
strictly pseudoconvex manifold with generating Levi
form defined by (1.1).
Then for every $\epsilon>0$ there exists $\delta>0$
such that for every
$\lambda\in\C^k$,
$c\in\R^k$,
$w_0\in\C^n$,
$y_0\in\R^k$,
$v\in\C^n$
such that
$$
\sum\Re(\lambda_j\z+c_j)A_j > \epsilon(|\lambda|+|c|)I
$$
and
$|w_0|<\delta,
|y_0|<\delta,
|v|<\delta$
there exists a unique stationary disc
$\z\mapsto f(\z)=(z(\z),w(\z))$
such that
$w(1)=w_0$,
$w'(1)=v$,
$y(1)=y_0$
that admits a lift $f^*$ such that
$f^*|_{b\Delta}=\Re(\lambda\zeta+c)G\partial\rho$
(where $\lambda$ and $c$ are handled as row-vectors,
and $G$ is a $k\times k$ matrix function
on $b\Delta$ close to the identity matrix uniquely
determined by $f$, see [T1]).
The pair $(f,f^*)$ depends smoothly (resp. analytically)
on $\zeta\in\bar\Delta$ and all the parameters
$\lambda, c, w_0, y_0, v$.
The pair $(f,f^*)$ is extremal in a suitable coordinate
system depending on $\epsilon$ only.
}
\medskip

Let $M$ be a generic manifold in $\C^N$ defined by (1.1).
Let $Q$ be the quadratic manifold obtained from (1.1)
by dropping the big `O' terms. We call $M$ {\it defective}
at 0 if all stationary discs for $Q$ are defective.
(That is every stationary disc, which possibly has a
lift with a pole at zero, also has another lift without
the pole. The authors do not know whether this situation
actually can occur.)
This definition is equivalent to the one given in [T2].
If the Levi form of $M$ splits into scalar forms, then
$M$ is not defective. Then for fixed $\epsilon$ and
sufficiently small $\delta$, all stationary discs provided by
Theorem 2.1 are not defective (see [T1], Proposition 6.8. or
[T2], Proposition 8.4).

Define $\L f={d\over d\theta}\big|_{\theta=0}f(e^{i\theta})$.
Note that if $f$ is holomorphic at $1\in\C$, then
$\L f=Jf'(1)$. Let $\zeta_0\in b\Delta$, $\zeta_0\ne1$.
Let $\E$ denote the set of all extremal pairs $(f,f^*)$
obtained by Theorem 2.1 such that $f$ is not defective.
If $M$ is a strictly convex hypersurface, then $\E$ stands for
the set of all extremal pairs, in which case $\E$ is a smooth
manifold by Lempert's [L] theory.
Consider the following evaluation maps:
$$
\eqalign{
&\F: \E\ni (f,f^*) \mapsto (f(1), f^*(1), \L f, \L f^*)\in TN^*(M), \cr
&\G: \E\ni (f,f^*) \mapsto (f(1), f^*(1), f(\zeta_0), f^*(\zeta_0)) \in
N^*(M)\times N^*(M).
}
$$

{\bf Proposition 2.2.} {\sl
The maps $\F$ and $\G$ are injective.}
\medskip

For the map $\F$ the proposition is proved in [T1], Proposition 3.9.
The proof for $\G$ is similar. We also need the following
stronger version.
\medskip

{\bf Proposition 2.2'.} {\sl
The maps $\F$ and $\G$ are diffeomorphisms onto their images.}
\medskip
{\it Proof.}
The source and target spaces of both $\F$ and $\G$ have the
same dimension $4N$. Hence it suffices to show that
$\F$ and $\G$ are immersions.
By an infinitesimal perturbation $(\dot f, \dot f^*)$
of an extremal pair $(f,f^*)$, we mean an element of the
tangent space to the finite dimensional manifold $\E$ at $(f,f^*)$.
To show that $\F$ is an immersion, we need to show that
$\dot f(1)=0$, $\dot f^*(1)=0$, $\L f=0$, and $\L f^*=0$
imply
$\dot f=0$ and $\dot f^*=0$.

We realize $(\dot f, \dot f^*)={d\over dt}\big|_{t=0}(f_t,f^*_t)$,
where $(f_t,f^*_t)$ is a smooth one parameter family of
extremal pairs with $(f_0,f^*_0)=(f,f^*)$.
For small $t$ all the pairs are close to $(f,f^*)$, hence
we can choose $\epsilon$ in (2.1) the same for all small $t$.
By (2.1) we have on $b\Delta$:
$$
\Re\<f^*_0,f_t-f_0\> \ge \epsilon |f_t-f_0|^2,\qquad
\Re\<f^*_t,f_0-f_t\> \ge \epsilon |f_0-f_t|^2.
$$
Adding the two inequalities yields
$$
\Re\<f^*_t-f^*_0,f_t-f_0\> \le -2\epsilon |f_t-f_0|^2.
$$
Dividing by $t^2$ and letting $t\to 0$ yields
$$
\Re\<\dot f^*,\dot f\> \le -2\epsilon |\dot f|^2
$$
for $\zeta\in b\Delta$.
The hypotheses imply
$\dot f=O(|\zeta-1|^2)$, $\dot f^*=O(|\zeta-1|^2)$.
Then
$$
\Re\int_0^{2\pi}{\<\dot f^*,\dot f\>\,d\theta\over |\zeta-1|^4} \le
-2\epsilon \int_0^{2\pi}{|\dot f|^2\,d\theta\over |\zeta-1|^4},
$$
where $\zeta=e^{i\theta}$.
Note for $|\zeta|=1$ we have
$d\zeta=i\zeta\,d\theta$ and
$\zeta|\zeta-1|^2=-(\zeta-1)^2$.
Then
$$
\int_0^{2\pi}{\<\dot f^*,\dot f\>\,d\theta\over |\zeta-1|^4}
= -i\int_{b\Delta}\left\<{\zeta \dot f^*\over(\zeta-1)^2},
{\dot f \over(\zeta-1)^2}\right\>\,d\zeta=0
$$
since the integrand is holomorphic in $\Delta$.
Hence
$$
\int_0^{2\pi}{|\dot f|^2\,d\theta\over |\zeta-1|^4}=0
$$
and $\dot f=0$.
Since $\dot f=0$, then  $\dot f^*|_{b\Delta}$ is tangent to the fibers
of $N^*(M)$ and gives rise to a lift of $f$.
Since
$\dot f^*=O(|\zeta-1|^2)$, then $\tilde f=\zeta(\zeta-1)^{-2}\dot f^*$
is a lift of $f$ without a pole at zero. Since $f$ is not defective,
then $\tilde f=0$, whence $\dot f^*=0$, and $\F$ is an immersion.
The proof that $\G$ is an immersion is similar.
It uses the identity
$\zeta\zeta_0|\zeta-\zeta_0|^2=-(\zeta-\zeta_0)^2$
for $|\zeta|=|\zeta_0|=1$.
The proof is complete.
\medskip

Define $T^+N^*(M)\subset TN^*(M)$.
We put
$\xi\in T^+_{(p,\phi)}N^*(M)$ if $\phi\in\Gamma^*_p$ and
$\pi_*\xi\in H_p$, where $\pi: T^*(\C^N)\to\C^N$
is the natural projection, and the Levi cones $\Gamma^*_p$ and $H_p$
are defined in Section 1.
\medskip

{\bf Proposition 2.3.} {\sl
Let $M$ be a strictly convex hypersurface in $\C^{n+1}$.
Then $\F(\E)=T^+N^*(M)$.}
\medskip
{\it Proof.}
The inclusion $\F(\E)\subset T^+N^*(M)$ follows by the Hopf lemma.
Indeed, let $M$ bound the domain $D$ defined by $\rho<0$, where
$\rho$ is a strictly convex function. Let $f$ be a nonconstant
analytic, not necessarily stationary disc attached to $M$,
and let $f(1)=p\in M$.
Then the nonconstant subharmonic function $\rho\circ f$ in $\Delta$
is zero on the boundary. By the Hopf lemma, $\<d\rho,f'(1)\>>0$.
This implies $-[f'(1)]\in\Gamma_p$ whence $\L f=Jf'(1)\in H_p$.
If $(f,f^*)\in\E$, then $f^*(1)\in\Gamma^*_p$,
and the desired inclusion follows.

The surjectivity of $\F$ follows by a simple topological argument.
Fix $p\in M$. Put $\E_p=\{(f,f^*)\in\E: f(1)=p\}$. Then the set $\E_p$
is contractible because  $f$ is completely determined by
$f(0)\in D$ and $f(1)=p$, and for given $f$, the supporting lift
$f^*$ is unique up to a positive multiplicative constant (see [L]).

Given $(f,f^*)\in\E_p$, we make a substitution by an automorphism
of the unit disc
$\zeta={\tau-\tau_0\over 1-\bar\tau_0\tau}e^{i\theta_0}$ with fixed point 1.
Put
$$
g(\tau)=f(\zeta),\qquad
g^*(\tau)=f^*(\zeta){(\tau-\tau_0)(1-\bar\tau_0\tau)\over\tau |1-\tau_0|^2},
$$
where we choose the factor so that $g^*$ has a pole at zero and
$g^*(1)=f^*(1)$.
Then $(g,g^*)\in\E_p$, and one can further check that
$$
\L g=\alpha \L f, \qquad
\L g^*=\L f^* - \beta f^*(1),
\eqno{(2.2)}
$$
where $\alpha,\beta\in\R$, $\alpha+i\beta={1+\tau_0\over 1-\tau_0}$.
Since $\tau_0\in\Delta$ is arbitrary, then $\alpha>0$ and $\beta\in\R$
are arbitrary.

Consider the map
$\Phi: \E_p\ni(f,f^*)\mapsto {\L f\over|\L f|}\in S^+$,
where $S^+=S^{2n+1}\cap H_p$ is the unit hemisphere in $H_p$.
By (2.2), the preimages of the map $\Phi$ are contractible.
Since $\E_p$ is contractible, so is $\Phi(\E_p)$.
It suffices to show that $\Phi(\E_p)=S^+$.
We will show that $\Phi(\E_p)$ contains an arbitrary small perturbation of
the equator of the hemisphere $S^+$. Then $\Phi(\E_p)$ will have to be
all of $S^+$.

We introduce a coordinate system $(z=x+iy,w)\in \C\times\C^n$ so that $p=0$ and
$M$ has a local equation
$$
x=|w|^2+O(|y|^3+|w|^3).
$$
Then $T_p(M)$ is defined by $x=0$ and $H_p\subset T_p(M)$ is the half-space $y<0$.
The stationary disc $f$ constructed by Theorem 2.1 for
$\lambda=0$, $c=1$, $w_0=0$, $y_0=0$, and small $v\in\C^n$
has the following asymptotic expression (see [T1], Corollary 5.2):
$$
z(\zeta)=O(|v|^2), \qquad
w(\zeta)=(\zeta-1)v +O(|v|^2).
$$
Then
$$
{\L f\over|\L f|}=\left(0,{v\over|v|}\right)+O(|v|), \quad |v|=\epsilon
$$
for small $\epsilon$ describes a small perturbation of
the equator of the hemisphere $S^+$.
Hence $\Phi(\E_p)=S^+$, and the proof is complete.
\medskip

If $M$ is a product of strictly convex hypersurfaces,
then $N^*(M)$, $T^+N^*(M)$, $\E$, etc., are the products of the
corresponding objects for the components of the product.
Then we immediately derive
\medskip

{\bf Corollary 2.4.} {\sl
Let $M$ be a product of strictly convex hypersurfaces.
Then $\F(\E)=T^+N^*(M)$.
}
\medskip

\beginsection  3. The main result

{\bf Theorem 3.1.} {\sl
Let $M_1$ be a real analytic strictly pseudoconvex generic manifold and
let $M_2$ be a product of several real analytic strictly convex hypersurfaces.
Then every biholomorphic map taking an open set in $M_1$ to $M_2$
continues along any path in $M_1$ as a locally biholomorphic map.
}
\medskip
{\it Remark.}
We require that $M_2$ be a product because we use Corollary 2.4
in the proof.
It would be interesting to find out for what manifolds the
conclusion of Corollary 2.4 is valid.

{\it Proof.}
The main idea of the proof is that a biholomorphism preserves
extremal pairs, therefore it extends along the extremal discs.

Let $F$ be a biholomorphic map defined at $p_1\in M_1$,  such that
$F(U)\subset M_2$ for some open set $U\subset M_1$.
The map $F$ lifts to the cotangent bundle $T^*(\C^N)$ in the usual way.
With some abuse of notation, we use the same letter $F$ for the lifted map.
We choose a coordinate system in which $p_1=0$ and $M_1$ is
given by (1.1).
Since $M_2$ is a product, then the Levi form of $M_2$ splits into
scalar positive definite forms.
Since the biholomorphic map $F$ preserves the Levi forms,
then the Levi form of $M_1$ at $p_1=0$ also splits into scalar positive
definite forms, and after a linear change of coordinates, the equation
of $M_1$ takes the form
$$
x_j=h_j(y,w)=|w_j|^2+O(|y|^3+|w|^3),\qquad
w_j\in\C^{n_j},\quad
n_1+\cdots+n_k=n.
\eqno{(3.1)}
$$
We note that the size of the coordinate chart for which (3.1)
holds is independent of the map $F$.
Indeed, once we know that such $F$ exists, then by Proposition 1.1
the Levi form of (the component of) $M_1$ splits onto scalar
forms. Then when extending $F$ along a path, we can always
restrict to finitely many coordinate charts by the compactness
argument.

Although Corollary 2.4 generally fails for $M_1$,
there are many extremal pairs $(f,f^*)$
such that $\F(f,f^*)\in T^+N^*(M_1)$.
Indeed, let
$(f,f^*)$ be the extremal pair constructed by Theorem 2.1 for
$\lambda=0$, $c_j=1$, $w_0=0$, $y_0=0$, and small $v\in\C^n$.
Then the components of $f$ admit
the following asymptotic expression (see [T1], Corollary 5.2):
$$
z(\zeta)=O(|v|^2), \qquad
w(\zeta)=(\zeta-1)v +O(|v|^2).
\eqno{(3.2)}
$$
Furthermore, plugging (3.2) in (3.1) and using the identity
$|\zeta-1|^2=-2\Re(\zeta-1)$ for $|\zeta|=1$, we obtain
$$
z_j(\zeta)=-2(\zeta-1)|v_j|^2 +O(|v|^3), \qquad
\L z_j=-2i|v_j|^2 +O(|v|^3).
$$
Note that the Levi cone $H_0$ of $M_1$ is defined by $x=0$, $y_j<0$.
Thus, if all $|v_j|$ are small and comparable, then
$\L f\in H_0$ and $\F(f,f^*)\in T^+N^*(M_1)$.
The same is true for all extremal pairs constructed with values
of parameters $\lambda$ and $c$ close to the above.

Consider all extremal pairs
$(f_1,f^*_1)$ for $M_1$ with fixed $f_1(1)=p_1=0$ and $f^*_1(1)$
such that $\F(f_1,f^*_1)\in T^+N^*(M_1)$.
Denote the set of such pairs by $\E_1$.

Put $\xi=\F(f_1,f^*_1)$. Since $F$ preserves the Levi forms,
then $F_*\xi\in T^+N^*(M_2)$.
By Corollary 2.4 there exists a unique extremal pair
$(f_2,f^*_2)$ for $M_2$ such that $\F(f_2,f^*_2)=F_*\xi$.
Fix $\zeta_0\in b\Delta$, $\zeta_0\ne1$.
We define
$\tilde F((f_1,f^*_1)(\zeta_0))=(f_2,f^*_2)(\zeta_0)$.
By Proposition 2.2', the map $\tilde F$ is a diffeomorphism
on the set $\{ (f_1,f^*_1)(\zeta_0): (f_1,f^*_1)\in\E_1 \}$.
Since all the objects are real analytic, then
$\tilde F$ is real analytic on an open set in $N^*(M)$.
Note that as $v\to 0$,
the pair $(f_1,f^*_1)$ shrinks into a point.
This implies that that the map $\tilde F$ agrees with
$F$ on an open set in $N^*(M)$ because $F$ preserves extremal pairs.
The extension preserves the fibers of $N^*(M)$ because
so does $F$. Hence, $\tilde F$ defines a real analytic
diffeomorphism on the set
$\{ (f_1(\zeta_0): (f_1,f^*_1)\in\E_1 \}\subset M$.
By varying $\zeta_0\in b\Delta$, we extend $\tilde F$ as a real analytic
diffeomorphism on the set
$V=\cup\{ (f_1(b\Delta\setminus\{1\}): (f_1,f^*_1)\in\E_1 \}\subset M$.
Since $\tilde F$ is real analytic and satisfies the tangential
Cauchy-Riemann equations on an open set in $M$, then it is CR
on the whole set $V\subset M$ where it is defined. Then by
real analyticity it further extends to a biholomorphic map in
a neighborhood of $V$ in $\C^N$.

Thus we conclude that $F$ extends as a biholomorphic map
along the boundaries of the extremal discs $f_1$.
By Proposition 2.2' (see also [T1], Corollary 5.6)
the directions of the boundary curves of the discs
$f_1$ span the tangent space $T_{p_1}(M_1)$. Then it follows
that all points within the same connected component can
be reached by moving along the boundaries of such discs,
and the theorem follows.

\beginsection References

\frenchspacing

\item{[A]}
H. Alexander, Holomorphic mappings from the ball and polydisc,
{\jour Math. Ann. \vol 209} (1974),
249--256.

\item{[BER]}
M. S. Baouendi, P. Ebenfelt, and L. Rothschild,
{\sl
Real Submanifolds in Complex Space and Their Mappings},
Princeton Mathematical Series, vol. 47,
Princeton University Press, Princeton, NJ, 1999.

\item{[EKV]}
V. V. Ezhov, N. G. Kruzhilin, and A. G. Vitushkin,
Extension of local mappings of pseudoconvex surfaces (Russian),
{\sl Dokl. Akad. Nauk SSSR \vol 270} (1983), 271--274.

\item{[HS]}
C. D. Hill and R. Shafikov,
Holomorphic correspondences between CR manifolds,
{\sl Indiana Univ. Math. J. \vol 54} (2005), 417--441.

\item{[L]}
L. Lempert,
La m\'etrique de Kobayashi et la repr\'esentation des domaines
sur la boule,
{\sl Bull. Soc. Math. France \bf 109} (1981), 427--474.

\item{[P]}
S. Pinchuk,
Holomorphic mappings of real-analytic hypersurfaces,
{\jour Mat. Sb. (N.S.) \vol 105} (147) (1978),
574--593 (Russian). English translation:
{\jour Math. USSR-Sb. \vol 34} (1978), 503--519.

\item{[Sc]}
A. Scalari,
Extremal discs and CR geometry, Ph.D. thesis,
University of Illinois at Urbana-Champaign, 2001.

\item{[T1]}
A. Tumanov,
Extremal discs and the regularity of CR mappings in higher
codimension, {\sl Amer. J. Math. \bf 123} (2001), 445--473.

\item{[T2]}
\bysame
Extremal discs and the geometry of CR manifolds,
In: {\sl Real Methods in Complex and CR Geometry},
CIME Session, Martina Franca, 2002
({\sl Lect. Notes in Math. \vol 1848}, 191--212)
Springer 2004.

\item{[W]}
S. M. Webster,
On the mapping problem for algebraic real hypersurfaces,
{\jour Invent. Math. \vol 43}
(1977), 53--68.

\bigskip\bigskip

{\auth Alberto Scalari},
Citigroup Centre, Canada Square, London E14 5LB, United Kingdom.
E-mail: {\it alberto.scalari@citigroup.com}
\bigskip
{\auth Alexander Tumanov},
Department of Mathematics, University of Illinois,
Urbana, IL 61801, U.S.A.
E-mail: {\it tumanov@uiuc.edu}.

\bye